\documentclass[11pt]{article}
\usepackage[table]{xcolor}
\usepackage{graphicx}
\usepackage{amsmath,amsfonts,amssymb}
\usepackage{graphicx}
\usepackage{multirow}
\usepackage{tikz,pgf}
\usepackage{url}
\usepackage{amsmath}
\usepackage{graphicx}
\usepackage{epsfig}
\usepackage{epstopdf}
\usepackage{color}
\usepackage{enumitem}
\def\disp{\displaystyle}

\def\tto{\;{\lower 1pt \hbox{$\rightarrow$}}\kern -10pt
\hbox{\raise 2pt \hbox{$\rightarrow$}}\;}

\def\Bar{\overline}
\def\ra{\rangle}
\def\la{\langle}

\def\epsilon{\varepsilon}

\def\h{\hfill\Box}
\def\R{\Bbb R}

\def\N{\Bbb N}
\def\ox{\bar{x}}

\def\oy{\bar{y}}
\def\oz{\bar{z}}

\def\gph{\mbox{\rm gph}}

\def\epi{\mbox{\rm epi}}

\def\dom{\mbox{\rm dom}}

\def\iri{\mbox{\rm iri}\,}

\def\h{\hfill\square}

\def\ph{\varphi}

\def\oR{\Bar{\R}}

\def\ph{\varphi}

\def\oR{\Bar{\R}}

\def\qri{\mbox{\rm qri}}
\setlist[enumerate,1]{itemsep=0.0ex,parsep=0.5ex,label={\rm(\alph*)},leftmargin=*, align=left}
\newcounter{lk}

\usepackage[colorlinks=true]{hyperref}

\begin{document}
\begin{center}
{\sc\bf Relationships between Polyhedral Convex Sets and Generalized Polyhedral Convex Sets}\footnote{This research is funded by the Vietnam Ministry of Education and Training under grant
	number B2022-CTT-06.}\\[1ex]

{Nguyen Ngoc Luan}\footnote{Department of Mathematics and Informatics, Hanoi National University of Education, 136 Xuan
Thuy, Hanoi, Vietnam (luannn@hnue.edu.vn).}, { Nguyen Mau   Nam}\footnote{Fariborz Maseeh Department of Mathematics and Statistics, Portland State University, Portland, OR
97207, USA (mnn3@pdx.edu). Research of this author was partly supported by the USA National Science Foundation under grant DMS-2136228.},
{ Nguyen Nang Thieu}\footnote{Institute of Mathematics, Vietnam Academy of Science and Technology,
Hanoi, Vietnam \& The State University of New York - SUNY, Korea (nnthieu@math.ac.vn).},
{ Nguyen Dong Yen}\footnote{Institute of Mathematics, Vietnam Academy of Science and Technology, 18 Hoang Quoc Viet,
Hanoi 10307 (ndyen@math.ac.vn).}\\[2ex]
\end{center}
\small{\bf Abstract.}  In this paper we  study some relationships between polyhedral convex sets (PCS) and generalized polyhedral convex sets (GPCS). In particular, we clarify by a counterexample that the necessary and sufficient conditions for the separation of a convex set and a PCS obtained by Kung Fu Ng and Wen Song in [Fenchel duality in finite-dimensional setting and its applications, {\em Nonlinear Anal.} {\bf 55} (2003), 845--858; Theorem~3.1] are no longer valid when considering GPCS instead of PCS. We also introduce and study the notions of generalized polyhedral set-valued mappings and optimal value functions generated by generalized polyhedral convex set-valued mappings along with their generalized differentiation calculus rules.   \\[1ex]

\newtheorem{Theorem}{Theorem}[section]
\newtheorem{Proposition}[Theorem]{Proposition}
\newtheorem{Remark}[Theorem]{Remark}
\newtheorem{Lemma}[Theorem]{Lemma}
\newtheorem{Corollary}[Theorem]{Corollary}
\newtheorem{Definition}[Theorem]{Definition}
\newtheorem{Example}[Theorem]{Example}
\renewcommand{\theequation}{\thesection.\arabic{equation}}
\normalsize

\section{Introduction}\label{Sect_1}
\setcounter{equation}{0}

Polyhedral convex sets (PCS) and related concepts have been study broadly in the framework of convex analysis in both finite dimensions and infinite dimensions.   Among many important results involving polyhedral convex sets,  necessary and sufficient conditions for the separation property of a convex set and polyhedral convex sets play a crucial role in developing generalized differentiation involving PCS with applications to optimization. An important result was established by Kung Fu Ng and Wen Song (see~\cite[Theorem~3.86]{mordukhovich_nam_2021},~\cite[Theorem~3.1]{ng-song}) providing necessary and sufficient conditions for separating a convex set and a polyhedral convex set in locally convex topological vector spaces. This is a generalization of a well-known result obtained by Rockafellar in finite dimensions; see~\cite[Theorem~20.2]{Rockafellar1970}. The result by Ng and Song was then used in~\cite{cuong-poly} to obtain comprehensive generalized differentiation calculus for nonsmooth functions and set-valued mappings in locally convex topological vector spaces.

Given the crucial role of PCS in convex analysis and applications, a new concept called {\em generalized polyhedral convex sets} (GPCS) have been introduced and studied recently in infinite dimensions; see~\cite{bs,luan_yen2020} and the references therein. In a series of recent papers, Luan, Yen and others have established the mathematical foundation for GPCS in locally convex topological vector spaces and provide important applications in many areas such as vector optimization, conic linear programming, numerical optimization, etc. These new developments also shed new light on many known results involving PCS. In particular, we refer the reader to the paper by Luan, Yao, and Yen~\cite{luan_yao_yen_NFAO2018} in which several  constructions such as sum of sets, sum of functions, directional derivative, infimal convolution, normal cone, subdifferential, conjugate function involving PCS and GPCS were thoroughly investigated.

The remarkable role of PCS and generalized PCS raises an important question asking for the clarification whether a certain result which holds for PCS also holds for GPCS or not.  One of the main goals of this paper  is to answer the question. In particular, we provide a counterexample showing that the aforementioned separation result by Ng and Song for PCS is no longer true for GPCS in general, and thus this counterexample somehow discourages the possibility for full generalizations of the results in~\cite{cuong-poly} to the case of GPCS. The second main goal of this paper is to study {\em polyhedral convex set-valued mappings} (PCSM) and {\em generalized polyhedral convex set-valued mappings} (GPCSM) and derive  calculus rules of generalized differentiation in  the case where one mapping involved is a PCSM, while the other mapping is a GPCSM. We also study generalized differentiation of optimal value functions generated by PCSM and GPCSM. When developing these generalized differentiation calculus rules, we recover a number of important results obtained by Luan, Yao, and Yen~\cite{luan_yao_yen_NFAO2018} by new proofs.

This paper is organized as follows. In Section~\ref{Sect_2} we study some relationships between PCS and GPCS in locally convex Hausdorff topological vector spaces. Section~\ref{Sect_3} is devoted to generalized differentiation for PCSM and GPCSM. Generalized differentiation of optimal value functions generated by PCSM and GPCSM is investigated in the final section.

\section{Relationships Between Polyhedral Convex Sets and Generalized Polyhedral Convex Sets}\label{Sect_2}
\setcounter{equation}{0}

Let $X$ be a locally convex Hausdorff topological vector space over the reals with its topological dual denoted by $X^*$. For simplicity of presentation, we assume that all spaces under consideration are locally convex Hausdorff topological vector spaces.  The \textit{cone} generated by a nonempty subset $\Omega$ of $X$ (resp., the \textit{closure} of $\Omega$) is denoted by ${\rm cone}(\Omega)$ (resp., $\overline{\Omega}$) . Thus, ${\rm cone}(\Omega)=\left\{tx\mid t\geq 0, x\in\Omega\right\}$.  In what follows, let $\oR:=(-\infty, \infty]$.

\begin{Definition}
\begin{enumerate}
\item A subset $P$ of $X$ is said to be a {\bf\em polyhedral convex set} (PCS) or a {\bf\em convex polyhedron} if there exist $x_1^*, \ldots, x^*_m\in X^*$ and $\alpha_1, \ldots, \alpha_m\in \R$ such that
    \begin{equation*}
    P=\big\{x\in X\; |\; \la x^*_i, x\ra\leq \alpha_i\; \mbox{\rm for all }i=1, \ldots, m\big\}.
    \end{equation*}
    \item A subset $Q$ of $X$ is said to be a {\bf\em generalized polyhedral convex set} (GPCS) or a {\bf\em generalized convex polyhedron} if there exist $x_1^*, \ldots, x^*_m\in X^*$,  $\alpha_1, \ldots, \alpha_m\in \R$, and a closed affine subspace $M$ of $X$  such that
        \begin{equation*}
    Q=\big\{x\in X\; |\; x\in M, \;\la x^*_i, x\ra\leq \alpha_i\; \mbox{\rm for all }i=1, \ldots, m\big\}.
    \end{equation*}
\end{enumerate}
\end{Definition}

It follows from the definitions that a GPCS can be represented as the intersection of a PCS and a closed affine subspace.

\begin{Proposition}\label{R1} Let $Q$ be a nonempty GPCS given by
 \begin{equation*}
    Q=\big\{x\in X\; |\; x\in M, \;\la x^*_i, x\ra\leq \alpha_i\; \mbox{\rm for all }i=1, \ldots, m\big\},
    \end{equation*}
    where $x_1^*, \ldots, x^*_m\in X^*$,  $\alpha_1, \ldots, \alpha_m\in \R$, and $M$ is a closed affine subspace. Then $Q$ is a PCS if and only if $\mbox{\rm codim}\, M<\infty$.
\end{Proposition}
{\bf Proof.} $\Longleftarrow$: Let $L:=M-M$ and observe that $L$ is the unique closed linear subspace parallel to $M$. Since $M$ has finite codimension, $\mbox{\rm dim}\, X/L=k$ for some positive integer~$k$. Consider the quotient mapping $\Phi\colon X\to X/L$. Choose $a\in M$ and get $[a]\in X/L$ with $\Phi(M)=[a]$. Suppose that $[b_1], \ldots, [b_k]$ form a basis for $X/L$. Then we have the representation
\begin{equation*}
[a]=\sum_{i=1}^k \beta_i [b_i],
\end{equation*}
where $\beta_1, \ldots, \beta_k\in \R$. For each $i\in\{1,\dots,k\}$, we consider the linear mapping $g_i\colon X/L\to\mathbb{R}$ with $g_i\big(\sum_{i=1}^{k}\mu_i[b_i]\big)=\mu_i$. By~\cite[Lemma~1.20 and Theorem~1.21]{rudinfunctional} (note that the same results and proofs are valid for linear mappings defined on $\mathbb R^n$), $g_i$ is continuous for all $i\in\{1,\dots,k\}$. Thus, the function $u_i^*:=g_i\circ\Phi$ is linear and continuous for all $i\in\{1,\dots,k\}$. Hence, $u_i^*\in X^*$ and $u_i^*\neq 0$ for all $i\in\{1,\dots,k\}$. We will now prove that
\begin{equation}\label{affine_polyhedral}
M=\big\{x\in X\mid u_i^*(x)=\beta_i,\ \, \forall i\in \{1,\dots,k\}\big\}.
\end{equation}
For all $x\in M$, we have $\Phi(x)=\Phi(a)=[a]=\sum_{i=1}^{k}\beta_i[b_i].$
Then, $u_i^*(x)=g_i(\Phi(x))=\beta_i$ for all $i\in\{1,\dots,k\}$. Therefore, $M\subset\{x\in X\mid u_i^*(x)=\beta_i,\ \, \forall i\in \{1,\dots,k\}\}.$ To prove the reverse inclusion, taking any $x$ in the right-hand side of~\eqref{affine_polyhedral}, one has $u_i^*(x)=\beta_i$ for all $i\in\{1,\dots,k\}$. Thus, $g_i(\Phi(x))=\beta_i$ for all $i\in\{1,\dots,k\}$, which means that $$\Phi(x)=\sum_{i=1}^{k}\beta_i[b_i]=\Phi(a).$$ Thus,  $\Phi(x-a)=[0]$. Hence, $x-a\in L$. This implies that $x\in M$. We have thus proved that the equality~\eqref{affine_polyhedral} is valid. Now, $Q$ can be represented as
\begin{equation*}
    Q=\big\{x\in X\; |\; \la u^*_i,x\ra=\beta_i\; \mbox{\rm for all }i=1, \ldots, k, \;\la x^*_i, x\ra\leq \alpha_i\; \mbox{\rm for all }i=1, \ldots, m\big\}.
    \end{equation*}
    Therefore, $Q$ is a PCS. \\[1ex]
    $\Longrightarrow$: Suppose that $Q$ is a PCS. Then there exist $z^*_1, \ldots, z^*_p\in X^*$ and $\gamma_1,\ldots, \gamma_p\in \R$ such that
    \begin{equation*}
    Q=\big\{x\in X\; |\; \la z^*_i, x\ra\leq \gamma_i\; \mbox{\rm for all }i=1, \ldots, p\big\}.
    \end{equation*}
    Choosing $x_0\in Q$ gives $x_0\in M$ and $\la z^*_i, x_0\ra\leq \gamma_i$ for all $i=1, \ldots, p$. Let $c_i:=\gamma_i-\la z^*_i, x_0\ra$ for $i=1, \ldots, p$.  Then $c_i\geq 0$ for all $i=1, \ldots, p$ and we have
    \begin{equation*}
    Q-x_0=\big\{y\in X\; |\; \la z^*_i, y\ra\leq c_i\big\}\subset L:=M-x_0.
    \end{equation*}
    Let $L_1:=\bigcap_{i=1}^p \mbox{\rm ker}\, z^*_i$. Then $L_1\subset Q-x_0\subset L$. It follows that
    \begin{equation*}
    \mbox{\rm codim}\, L\, \leq\, \mbox{\rm codim}\, L_1<\infty,
 \end{equation*}
 which completes the proof. $\h$

 The next corollary is a direct consequence of Proposition~\ref{R1}.

\begin{Corollary}\label{prop: affine_subspace_polyhedral}
	A closed linear subspace $M$  of $X$ is a PCS if and only if $M$ has finite codimension.
\end{Corollary}

Next, we present an example of a GPCS that is not a PCS.

\begin{Example}{\rm Let $X=\ell_2$ and let
		\begin{equation*}
			Q=\{x=(x_n)\; |\; x_{2k}=0\; \mbox{\rm for all }k\in \N\}.
		\end{equation*}
		Then $Q$ is a GPCS. In fact, it is a closed linear subspace of $X$. Since $\mbox{\rm codim}\, Q=\infty$,  the set $Q$ is not a PCS by Proposition~\ref{R1}.
}\end{Example}

\begin{Definition}\label{Poly Sep} Let $\Omega_1$ and $\Omega_2$ be a nonempty convex sets in $X$. We say that $\Omega_1$ and $\Omega_2$  can be separated by a closed hyperplane that does not contain $\Omega_2$ if there exist $x^*\in X^*$ and $\alpha\in \R$ such that
	\begin{equation*}
		\la x^*, x\ra\leq \alpha\leq \la x^*, y\ra\ \, \mbox{\rm whenever}\ x\in \Omega_1,\; y\in \Omega_2,
	\end{equation*}
	and there exists $\hat{y}\in \Omega_2$ such that $\alpha<\la x^*, \hat{y}\ra.$
\end{Definition}
In the setting of Definition~\ref{Poly Sep}, define
\begin{eqnarray*}
	\begin{array}{rcl}
		\mathcal{H}&=&\{x\in X^*\; |\; \la x^*, x\ra=\alpha\},\\
		\mathcal{H}_{+}&=&\{x\in X^*\;  |\; \la x^*, x\ra\geq\alpha\},\\
		\mathcal{H}_{-}&=&\{x\in X^*\;  |\; \la x^*, x\ra\leq\alpha\}.
	\end{array}
\end{eqnarray*}
Since $x^*$ is obviously nonzero, $\mathcal{H}$ is a closed hyperplane. We have
\begin{equation*}
	\Omega_1\subset\mathcal{H}_{-}, \;\Omega_2\subset \mathcal{H}_{+},\; \Omega_2\not\subset \mathcal{H}.
\end{equation*}

\begin{Definition} {\rm (See~\cite[Definition~2.168]{mordukhovich_nam_2021})} Let $\Omega$ be a convex subset of $X$.
	\begin{enumerate}
	\item The intrinsic relative interior of $\Omega$ is the set
	$${\rm iri}(\Omega) :=\left\{x\in \Omega \mid {\rm cone}(\Omega - x)\ \, {\rm is\ a\ linear\ subspace\ of}\ \, X\right\}.$$
	\item The quasi-relative interior of $\Omega$ is the set
	$${\rm qri}(\Omega):=\left\{x\in \Omega \mid \overline{\rm cone}(\Omega - x)\ \, {\rm is\ a\ linear\ subspace\ of}\ \, X\right\}.$$
\item We say that $\Omega$ is quasi-regular if ${\rm qri}(\Omega) = {\rm iri}(\Omega)$.
	\end{enumerate}
\end{Definition}

In~\cite{luan_yao_yen_NFAO2018} and the references therein, several important results for PCS  have been generalized for GPCS. We present below a number of important results which hold for PCS but do not hold for GPCS. The first one is a convex separation theorem involving a PCS and a convex set in $X$; see~\cite[Theorem~3.86]{mordukhovich_nam_2021}.

\begin{Theorem}\label{thm: subspace_separation}
	Let $P$ be a nonempty PCS and let $\Omega$ be a nonempty convex set in $X$. Suppose that  ${\rm qri}(\Omega)\neq\emptyset$. Then $P$ and $\Omega$ can be separated by a closed hyperplane that does not contain $\Omega$ if and only if $P\cap {\rm qri}(\Omega)=\emptyset$.
\end{Theorem}
This result plays a crucial role in developing generalized differentiation for nonsmooth functions and set-valued mappings in the case where some functions and mappings involved are generated by PCS.

The following example will show that the conclusion of Theorem~\ref{thm: subspace_separation}  may not hold true when  instead of $P$ one takes a subspace $M$  of infinite codimension.

\begin{Example}\label{ex: separation_subspace_fail}{\rm
		 Let $\Omega_0=\left\{x\in \ell_2\mid x \;\text{has finitely many nonzero coordinates}\right\}$, $$y=\left(1,\dfrac{1}{2},\dfrac{1}{4}, \dfrac{1}{8},\dots \right),\ \; z=\left(1,\dfrac{1}{2},\dfrac{1}{3},\dots\right),$$ and $M=\left\{\mu z\mid \mu\in \mathbb{R}\right\}$. Clearly, $y\in\ell_2$ and $z\in\ell_2$. Put $\Omega=y+\Omega_0$. Then, $\Omega$ is an affine subset of $\ell_2$. Hence, $\Omega$ is convex and $\qri(\Omega)=\iri(\Omega)=\Omega$. We have $M\cap \Omega=\emptyset$. Indeed, if $M\cap \Omega\neq\emptyset$, then for some $u\in\Omega_0$, we have $y+u\in M$. Thus, there is $\mu\in\R$ and $\bar{k}\in\mathbb{N}$ such that $\dfrac{1}{2^k}=\dfrac{\mu}{k+1}$ for all $k\geq \bar{k}$, which is a contradiction. As $\Omega=\qri(\Omega)$, it follows that $M\cap\qri(\Omega)=\emptyset$.
		
		Next, we will prove that there is no closed hyperplane which does not contain $\Omega$ and separates $M$ and $\Omega$. Suppose on the contrary that there exists a closed  hyperplane $\mathcal{H}\subset X$ such that $\mathcal{H}$ does not contain $\Omega$ and separates $M$ and $\Omega$. Then there is $x^*\in \ell_2\setminus \{0\}$ and $\alpha\in\R$ such that $\mathcal{H}=\{x\in\ell_2\mid \langle x^*,x\rangle = \alpha\}$. Since $\mathcal{H}$ separates $\Omega$ and $M$ and $\mathcal{H}$ does not contain $\Omega$, we have
		\begin{equation}\label{separation_inq}
			\sup\limits_{x\in M} \langle x^*,x\rangle \leq\alpha\leq \inf\limits_{x\in \Omega}\langle x^*,x\rangle
		\end{equation}
		and there exists $w\in\Omega$ such that $\alpha< \langle x^*,w\rangle.$  As $\langle x^*,x\rangle \leq\alpha$ for all $x\in M$, one has $\langle x^*,\mu z\rangle \leq\alpha$ for all $\mu\in\R$. Thus, $\langle x^*, z\rangle =0$. Then, the relation~\eqref{separation_inq} is equivalent to
		$$0 \leq\alpha\leq \inf\limits_{x\in \Omega}\langle x^*,x\rangle.$$
		 So, we have
		$$0 \leq\alpha\leq \langle x^*,y\rangle+\langle x^*,u\rangle$$
		for all $u\in\Omega_0$. This implies
		$$-\langle x^*,y\rangle\leq\langle x^*,u\rangle$$
		for all $u\in\Omega_0$. Since $\Omega_0$ is dense in $\ell_2$, the latter  property yields
		$$-\langle x^*,y\rangle\leq\langle x^*,v\rangle$$
		for all $v\in\ell_2$. This means that $x^*=0$, which contradicts the choice of $x^*$. Hence, there is no closed hyperplane not containing $\Omega$ which separates $M$ and $\Omega$.
	}
\end{Example}

\begin{Remark}{\rm
		In Example~\ref{ex: separation_subspace_fail}, the set $\Omega$ is not closed.  However,
		the assertion of Theorem~\ref{thm: subspace_separation}  may still be false when a subspace $M$ of infinite codimension plays the role of $P$  and $\Omega$ is a closed set in $X$. Indeed, by~\cite[Remark~2.12]{luan_yao_yen_NFAO2018}, there exists a locally convex topological vector space $X$ with two closed linear subspaces $L$ and $M$ such that $\overline{L+M}=X$ but $L+M\neq X$. Taking any $a\in X\setminus(L+M)$ and setting $\Omega=a+ L$, one  sees that $\Omega$ is a closed affine set. Hence, $\iri(\Omega)=\qri(\Omega)=\Omega$. First, we will show that $ M\cap\Omega=\emptyset$. Suppose on the contrary that there is some $u\in  M\cap\Omega$. Then, $u\in  M$ and $u=a+v$ for some $v\in L$. Therefore, $a=u-v\in L+M$, which is a contradiction. Thus, $ M\cap\Omega=\emptyset$. Next, we will prove that $ M$ and $\Omega$ cannot be separated by any hyperplane. Suppose on the contrary that there exist a nonzero linear  functional $x^*\in X^*$ and a real number $\alpha$ such that
		\begin{equation}\label{separation_inq2}
			\sup\limits_{x\in  M} \langle x^*,x\rangle \leq \alpha\leq \inf\limits_{y\in \Omega} \langle x^*,y\rangle.
		\end{equation}
		If there is some $\bar x\in   M$ such that $\langle x^*, \bar x\rangle =\beta\neq0$, then by taking $t=(\alpha+1)/\beta$, we have $t\bar x\in  M$ and $\langle x^*,t \bar x\rangle =\alpha+1> \alpha,$
		which  contradicts the fact that $\sup\limits_{x\in  M}\, \langle x^*,x\rangle \leq \alpha$. Hence $\langle x^*,x\rangle =0$ for all $x\in  M$. Thus,~\eqref{separation_inq2} yields
		$$0\leq \langle x^*,a\rangle +\inf\limits_{z\in  L}\, \langle x^*, z\rangle.$$
		Since $\langle x^*,a\rangle$ is fixed, the latter implies that  $\langle x^*,  z\rangle =0$ for all $ z\in  L$. Therefore, we have shown that $\langle x^*,x+z\rangle =0$ for all $x\in  M$ and $ z\in L$. Recalling that $\overline{L+M}=X$, we  can infer that $\langle x^*,  u\rangle =0$ for all $ u\in X$. This  contradicts the choice of $x^*$. We have thus proved that $ M\cap\qri(\Omega)=\emptyset$ and that $ M$ and $\Omega$ cannot be separated by any hyperplane.}
\end{Remark}

\section{Generalized Differentiation for Convex Polyhedral Set-Valued Mappings and  Generalized Convex Polyhedral Set-Valued Mappings}\label{Sect_3}
\setcounter{equation}{0}

The counterexample in Example~\ref{ex: separation_subspace_fail} shows that  analogues of the generalized differentiation results involving PCS in~\cite{mordukhovich_nam_2021}  may not hold for GPCS. This section   establishes some positive results.

\begin{Lemma}\label{LM1} Let $P$ be a   polyhedral convex set and let $M$ be a closed affine subspace in $X$. Then
	\begin{equation*}
		N(\ox; P\cap M)=N(\ox; P)+N(\ox; M)=N(\ox; P)+L^\perp\ \; \mbox{\rm for all }\; \ox\in P\cap M,
	\end{equation*}
	where $L$ is  the linear subspace parallel to $M$.
\end{Lemma}
{\bf Proof.} Fix any $\ox\in P\cap M$. Then $\ox\in P\cap \mbox{\rm qri}(M)=P\cap M$, so $P\cap \mbox{\rm qri}(M)\neq\emptyset$. It is follows from~\cite[Theorem~3.87]{mordukhovich_nam_2021} that
\begin{equation*}
	N(\ox; P\cap M)=N(\ox; P)+N(\ox; M).
\end{equation*}
Since $N(\ox; M)=L^\perp$, this completes the proof. $\h$

\begin{Lemma}\label{LM2} Let $P_1$ and $P_2$ be two convex polyhedra in $X$. Then
	\begin{equation*}
		N(\ox; P_1\cap P_2)=N(\ox; P_1)+N(\ox; P_2)\ \; \mbox{\rm for all }\; \ox\in P_1\cap P_2.
	\end{equation*}
\end{Lemma}
{\bf Proof.} This is obvious because  if
$$P=\left\{x\in X\; |\; \la x^*_i, x\ra\leq \alpha_i\ \; \mbox{\rm for all }\; i=1, \ldots, m\right\},$$
then $N(\ox; P)=\mbox{\rm cone}\{x^*_i\; |\; i\in I(\ox)\}$, where
$I(\ox):=\{i \; |\; i=1, \ldots, m,\ \la x^*_i, \ox\ra=\alpha_i\}$. $\h$

Next, we present  a new proof for the important result  obtained by Luan, Yao, and Yen in~\cite[Theorem~4.10]{luan_yao_yen_NFAO2018}.

\begin{Theorem}\label{P intersection rule} Let $P$ be a  PCS and let $Q$ be a GPCS. Then
	\begin{equation*}
		N(\ox; P\cap Q)=N(\ox; P)+N(\ox; Q)\ \; \mbox{\rm for all }\; \ox\in P\cap Q.
	\end{equation*}
\end{Theorem}
{\bf Proof.}  Fix a point $\ox\in P\cap Q$. Since $Q$ is a GPCS, we have the representation
\begin{equation*}
	Q=P_1\cap M,
\end{equation*}
where  $P_1$ is a PCS and $M$ is a closed affine subspace. Then by Lemma~\ref{LM1} and Lemma~\ref{LM2} we have
\begin{align*}
	N(\ox; P\cap Q)&=N(\ox; (P\cap P_1)\cap M)=N(\ox; P\cap P_1)+N(\ox; M)\\
	&=N(\ox; P)+N(\ox; P_1)+N(\ox; M)=N(\ox; P)+N(\ox; P_1\cap M)\\
&=N(\ox; P)+N(\ox; Q).
\end{align*}
This completes the proof. $\h$

\medskip
In what follows, let $X$ and $Y$ be locally convex Hausdorff topological vector spaces over the reals. For a set-valued mapping $F\colon X\tto Y$, one defines the \textit{graph} and the \textit{effective domain} of~$F$ respectively by $$\gph(F):=\big\{(x,y)\in X\times Y\mid y\in F(x)\big\}$$ and $$\dom(F):=\big\{x\in X\mid F(x)\neq\emptyset\big\}.$$

\begin{Definition}   Let $F\colon X\tto Y$ be a set-valued mapping.
	\begin{enumerate}
		\item $F$ is said to be {\bf\em convex} if $\gph(F)$ is a convex set in $X\times Y$.
		\item $F$ is said to be {\bf\em polyhedral convex} if $\gph(F)$ is a PCS in $X\times Y$.
		\item  $F$ is said to be {\bf\em generalized  polyhedral convex} if $\gph(F)$ is a GPCS in $X\times Y$.
	\end{enumerate}
\end{Definition}

\begin{Definition}\label{coder}  Let $F\colon X\tto Y$ be a convex set-valued mapping and let $(\ox,\oy)\in \gph(F)$. The {\bf\em coderivative} of $F$ at $(\ox,\oy)$ is the set-valued mapping $D^*F(\ox,\oy)\colon Y^*\tto X^*$ with the values
	\begin{equation}\label{cod}
		D^*F(\ox,\oy)(v^*):=\big\{u^*\in X^*\;\big|\;(u^*,-v^*)\in N\big((\ox,\oy);\gph(F)\big)\big\},\ \, v^*\in Y^*.
	\end{equation}
\end{Definition}

\begin{Example}\label{indicator mapping} {\rm Given a subset $\Theta$ of $X$, define $\Delta_\Theta\colon X\tto Y$ by
		\begin{equation*}
			\Delta_\Theta(x):=  \begin{cases} 0&\mbox{\rm if }x\in \Theta,\\
				\emptyset&\mbox{\rm if }x\notin \Theta.
			\end{cases}
		\end{equation*}
		Then $\gph(\Delta_\Theta)=\Theta\times \{0\}$. Suppose that $\Theta$ is a convex set  and $\ox\in \Theta$. We have $N((\ox, 0); \gph(\Delta_\Theta))=N(\ox; \Theta)\times Y$ and hence
		\begin{equation}\label{coderivative_Delta}
			D^*\Delta_\Theta(\ox, 0)(v^*)=N(\ox; \Theta)\ \; \mbox{\rm for all }\; v^*\in Y^*.
		\end{equation}
	}
\end{Example}

For any convex set-valued mappings $F_1,F_2\colon X\tto Y$, it follows from the definition that their  sum, which is defined by setting $(F_1+F_2)(x)=F_1(x)+F_2(x)$ for all $x\in X$, is a convex set-valued mapping with $\mbox{\rm dom}(F_1+F_2)=\dom(F_1)\cap\dom(F_2)$. Our first calculus result concerns representing the coderivative of  $F_1+F_2$ at a given point $(\ox,\oy)\in\mbox{\rm gph}(F_1+F_2)$. To formulate this result, consider the nonempty set
\begin{equation}\label{S}
	S(\ox,\oy):=\big\{(\oy_1,\oy_2)\in Y\times Y\;\big|\;\oy=\oy_1+\oy_2,\;\oy_i\in F_i(\ox),\;i=1,2\big\}.
\end{equation}
The following theorem gives us the coderivative sum rule for  polyhedral convex set-valued mappings and generalized polyhedral convex  set-valued mappings.

\begin{Theorem}\label{CSR} Let $F_1\colon X\tto Y$ be a polyhedral convex set-valued mapping and let $F_2\colon X\tto Y$ be a generalized polyhedral convex set-valued mapping. Then  the equality
	\begin{equation*}
		D^*(F_1+F_2)(\ox, \oy)(v^*)=D^*F_1(\ox, \oy_1)(v^*)+D^*F_2(\ox, \oy_2)(v^*)
	\end{equation*}
holds for every $v^*\in Y^*$ whenever $(\oy_1, \oy_2)\in S(\ox, \oy)$, where $S$ is defined in~\eqref{S}.
\end{Theorem}
{\bf Proof.}   Let $(\oy_1, \oy_2)\in S(\ox, \oy)$ and $v^*\in Y^*$ be given arbitrarily. Fix any \begin{equation}\label{simple_inclusion} u^*\in D^*(F_1+F_2)(\ox,\oy)(v^*).\end{equation}
Then the inclusion $(u^*,-v^*)\in N((\ox,\oy);\mbox{\rm gph}(F_1+F_2))$  is valid. Consider the convex sets
\begin{eqnarray*}
	\begin{array}{ll}
		&\Omega_1:=\big\{(x,y_1,y_2)\in X\times Y\times Y\;\big|\;y_1\in F_1(x)\big\},\\
		&\Omega_2:=\big\{(x,y_1,y_2)\in X\times Y\times Y\;\big|\;y_2\in F_2(x)\big\}
	\end{array}
\end{eqnarray*}
and deduce from the normal cone definition that
\begin{equation}\label{inclusion_sum_rule}
	(u^*,-v^*,-v^*)\in N((\ox,\oy_1,\oy_2);\Omega_1\cap\Omega_2).
\end{equation}
Observe that $\Omega_1=\big(\gph F_1\big)\times Y$ and thus it is a PCS in $X\times Y\times Y$  by the assumption made on $F_1$. Similarly, $\Omega_2$ is a GPCS in $X\times Y\times Y$  by the condition imposed on $F_2$.
Then we can employ Theorem~\ref{P intersection rule} and get
\begin{equation*}
	(u^*,-v^*,-v^*)\in N((\ox,\oy_1,\oy_2);\Omega_1\cap\Omega_2)=N((\ox,\oy_1,\oy_2);\Omega_1)+N((\ox,\oy_1,\oy_2);\Omega_2).
\end{equation*}
Therefore, we  obtain the relationships
\begin{equation*}
	(u^*,-v^*,-v^*)=(u_1^*,-v^*,0)+(u_2^*,0,-v^*),
\end{equation*}
where $(u_i^*,-v^*)\in N((\ox,\oy_i);\gph F_i)$ for $i=1, 2$. This implies by the coderivative definition that
\begin{equation*}
	u^*=u_1^*+u_2^*\in D^*F_1(\ox,\oy_1)(v^*)+D^*F_2(\ox,\oy_2)(v^*).
\end{equation*}
So, we have proved that
$$D^*(F_1+F_2)(\ox,\oy)(v^*)\subset D^*F_1(\ox,\oy_1)(v^*)+D^*F_2(\ox,\oy_2)(v^*).$$
To prove the reverse inclusion, take any $u_1^*D^*F_1(\ox,\oy_1)(v^*)$ and $u_2^*\in D^*F_2(\ox,\oy_2)(v^*)$. Then $$\big\langle (u^*_i,-v^*), (x-\bar x, y_i-\oy_i)\big\rangle\leq 0$$ for every $x\in X$ and $y_i\in F_i(x)$ with $i=1,2$. It follows that
$$\big\langle (u^*_1,-v^*,0), (x-\bar x, y_1-\oy_1, y_2-\oy_2)\big\rangle\leq 0$$ and
$$\big\langle (u^*_1,0,-v^*), (x-\bar x, y_1-\oy_1, y_2-\oy_2)\big\rangle\leq 0.$$ Adding these inequalities side-by-side yields
$$\big\langle (u^*,-v^*,-v^*), (x-\bar x, y_1-\oy_1, y_2-\oy_2)\big\rangle\leq 0$$ for every $x\in X$ and $y_i\in F_i(x)$ with $i=1,2$, where $u^*:=u^*_1+u^*_2$. Hence one gets the inclusion~\eqref{inclusion_sum_rule}, which clearly implies that $(u^*,-v^*)\in N((\ox,\oy);\mbox{\rm gph}(F_1+F_2))$. Hence~\eqref{simple_inclusion} is valid, and thus we have verified the claimed sum rule. $\h$

\begin{Definition} Let $f\colon X\to \oR  =(-\infty, \infty]$ be an extended-real-valued function. The epigraph of $f$ is the set
	$$\epi(f):=\big\{(x,\alpha)\in X\times\mathbb R\mid \alpha\geq f(x)\big\}.$$
	\begin{enumerate}
		\item We say that $f$ is  {\bf\em polyhedral convex} if $\epi(f)$ is a PCS in  $X\times \R$.
		\item  We say that $f$ is  {\bf\em generalized  polyhedral convex} if $\epi(f)$ is a GPCS in $X\times\R$.
	\end{enumerate}
\end{Definition}

The \textit{effective domain} of an extended-real-valued function $f\colon X\to \oR$ is the set $$\dom(f):=\big\{x\in X\mid f(x)<\infty\big\}.$$ If $f$ is convex, then the \textit{subdifferential} $\partial f(\bar x)$ of $f$ at $\bar x\in\dom(f)$ is defined by setting 	\begin{eqnarray*}\begin{array}{rcl}
	\partial f(\ox) & = &  \big\{x^*\in X^*\mid \langle x^*,x-\ox\rangle\leq f(x)-f(\ox)\ \, {\rm for\ all\ }\, x\in X\big\}\\
		& = & \big\{x^*\in X^*\mid (x^*,-1)\in N((\ox, f(\ox));\epi(f))\big\}.
	\end{array}
\end{eqnarray*}
Theorem~\ref{CSR} allows us to obtain the next subdifferential sum rule for polyhedral convex functions and generalized polyhedral convex functions.

\begin{Corollary}\label{ssr} Let $f_1, f_2\colon X\to \oR$ be two extended-real-valued functions. Suppose that $f_1$ is a polyhedral convex function and $f_2$ is a generalized polyhedral convex function. Then
	\begin{equation*}
		\partial (f_1+f_2)(\ox)=\partial f_1(\ox)+\partial f_2(\ox)\ \; {\rm for \ every\ }\; \ox\in \dom(f_1)\cap \dom(f_2).
	\end{equation*}
\end{Corollary}
{\bf Proof.} Fix any $\ox\in \dom(f_1)\cap \dom(f_2)$. Let $F_i(x):=[f_i(x), \infty)$  for all $x\in X$ and get $\gph(F_i)=\epi(f_i)$ for $i=1, 2$. Thus, $F_1$ is a polyhedral convex set-valued mapping and $F_2$ is a generalized polyhedral convex set-valued mapping. In addition,
\begin{equation*}
	D^*F_i(\ox, f_i(\ox))(1)=\partial f_i(\ox)\ \; \mbox{\rm for }\; i=1, 2.
\end{equation*}
Let $\oy:=f_1(\ox)+f_2(\ox)$. Then $S(\ox, \oy)=\{\big(f_1(\ox), f_2(\ox)\big)\}$, where $S(\ox, \oy)$ is defined in~\eqref{S}. Applying Theorem~\ref{CSR} gives
\begin{align*}
	\partial (f_1+f_2)(\ox)&=D^*(F_1+F_2)(\ox, \oy)(1)=D^*F_1(\ox, f_1(\ox))(1)+D^*F_2(\ox, f_2(\ox))(1)\\
	&=\partial f_1(\ox)+\partial f_2(\ox).
\end{align*}
This completes the proof. $\h$

Now we consider the  composition of two mappings $F\colon X\tto Y$ and $G\colon Y\tto Z$.
It follows from the definition that $G\circ F$ is convex provided that both $F$ and $G$ have this property. Given $\oz\in(G\circ F)(\ox)$, we consider the set
\begin{equation}\label{M}
	M(\ox,\oz):=F(\ox)\cap G^{-1}(\oz).
\end{equation}
The following theorem establishes the coderivative chain rule for set-valued mappings.

\begin{Theorem}\label{scr} Let $F\colon X\tto Y$ and $G\colon Y\tto Z$ be  set-valued mappings. Suppose that $F$ is a polyhedral set-valued mapping and $G$ is a generalized polyhedral set-valued mapping or vice versa.
	Then for any $(\ox,\oz)\in\mbox{\rm gph}(G\circ F)$ and $ w^*\in Z^*$ we have the coderivative chain rule
	\begin{equation}\label{chain}
		D^*(G\circ F)(\ox,\oz)(w^*)=D^*F(\ox,\oy)\circ D^*G(\oy,\oz)( w^*)
	\end{equation}
	whenever $\oy\in M(\ox,\oz)$.
\end{Theorem}
{\bf Proof.}  Picking $u^*\in D^*(G\circ F)(\ox,\oz)( w^*)$ and $\oy\in M(\ox,\oz)$ gives us the inclusion $$(u^*,- w^*)\in N((\ox,\oz);\mbox{\rm gph}(G\circ F)),$$ which means that
\begin{equation}\label{equation_1}
	\la u^*,x-\ox\ra-\la w^*,z-\oz\ra\le 0\ \;\mbox{\rm for all }\;(x,z)\in\mbox{\rm gph}(G\circ F).
\end{equation}
Define two convex subsets of $X\times Y\times Z$ by
$$
\Omega_1:=(\gph F)\times Z\;\mbox{ and }\;\Omega_2:=X\times(\gph G).
$$
We can directly deduce  from~\eqref{equation_1} and the definitions that
\begin{equation}\label{equation_2}
	(u^*,0,-w^*)\in N((\ox,\oy,\oz);\Omega_1\cap\Omega_2).
\end{equation}
Applying Theorem~\ref{P intersection rule}  together with the conditions made on $F_1$ and $F_2$ tells us that
\begin{equation}\label{equation_3}
	(u^*,0,-w^*)\in N((\ox,\oy,\oz);\Omega_1\cap\Omega_2)=N((\ox,\oy,\oz);\Omega_1)+N((\ox,\oy,\oz);\Omega_2),
\end{equation}
and thus there exists a vector  $v^*\in Y^*$ such that we have the representation
\begin{equation}\label{equation_3a}
	(u^*,0,-w^*)=(u^*, -v^*,0)+(0,v^*,-w^*)
\end{equation}
with $(u^*,-v^*)\in N((\ox,\oy);\gph F)$  and $(v^*,-w^*)\in N((\oy,\oz);\gph G)$.
This shows by the coderivative definition  in~\eqref{cod} that
\begin{equation}\label{equation_4}
	u^*\in D^*F(\ox,\oy)(v^*)\ \;\mbox{\rm and }\; v^*\in D^*G(\oy,\oz)(w^*),
\end{equation}
and so we  get the inclusion ``$\subset$" in \eqref{chain}.  The reverse inclusion can be proved as follows. Given any $u^*\in D^*F(\ox,\oy)\circ D^*G(\oy,\oz)(w^*)$, one can find some $v^*\in Y^*$ such that~\eqref{equation_4} holds. Then~\eqref{equation_3a} is fulfilled and, moreover, one has $(u^*, -v^*,0)\in N((\ox,\oy,\oz);\Omega_1)$ and $(0,v^*,-w^*)\in N((\ox,\oy,\oz);\Omega_2)$. Since the inclusion
$$N((\ox,\oy,\oz);\Omega_1)+N((\ox,\oy,\oz);\Omega_2) \subset   N((\ox,\oy,\oz);\Omega_1\cap\Omega_2)$$ is valid whenever $F$ and $G$ are merely convex set-valued mappings, one gets~\eqref{equation_2}, which implies~\eqref{equation_1}. Hence, $u^*\in D^*(G\circ F)(\ox,\oz)(w^*)$. The proof is complete. $\h$

\medskip
The next rule for computing subdifferentials of the composition of a polyhedral convex function and an affine mapping is a corollary of the preceding theorem.

\begin{Corollary} Let $B\colon X\to Y$ be an affine mapping given by
	\begin{equation*}
		B(x):=A(x)+b\ \; \mbox{\rm for }\; x\in X,
	\end{equation*}
	where $A\colon X\to Y$ is a continuous linear mapping and $b\in Y$. If $f\colon Y\to \oR$ is a polyhedral convex function, then
	\begin{equation}\label{subdifferential_of_composition}
		\partial (f\circ B)(\ox)= A^*\big(\partial f(\oy)\big)\ \; \mbox{\rm for all }\; \ox\in \dom(f\circ B),
	\end{equation}
	where $\oy:=B(\ox)$.
\end{Corollary}
{\bf Proof.} Let $F(x):=\{B(x)\}$ for $x\in X$ and let $G(y):=[f(y), \infty)$ for $y\in Y$. Then
\begin{equation}\label{composition_G_F}
	(G\circ F)(x)=[(f\circ B)(x), \infty)\ \; \mbox{\rm for all }\; x\in X.
\end{equation}
By our assumptions, $F:X\rightrightarrows Y$ is a generalized polyhedral set-valued mapping and $G:Y\rightrightarrows\mathbb R$ is a polyhedral set-valued mapping. Applying the coderivative chain rule from Theorem~\ref{scr}, we can get the desired result. Indeed, take any $\ox\in \dom(f\circ B)$ and put $\oy=B(\ox)$, $\oz=f(\oy)$. Then $(\ox,\oz)\in\mbox{\rm gph}(G\circ F)$ and $\oy\in M(\ox,\oz)=F(\ox)\cap G^{-1}(\oz)$. Therefore, by~\eqref{chain} we have
\begin{eqnarray*}\begin{array}{rcl}
			D^*(G\circ F)(\ox,\oz)(-1)  & = & D^*F(\ox,\oy)\circ D^*G(\oy,\oz)(-1)\\
		& = & D^*F(\ox,\oy)\left(\partial f(\oy)\big)\right)\\
		& = & A^*\big(\partial f(\oy)\big),
	\end{array}
\end{eqnarray*}
which together with~\eqref{composition_G_F} implies the equality in~\eqref{subdifferential_of_composition}.
$\h$

\medskip
Let $F\colon X\tto Y$ be a set-valued mapping and let $\Theta\subset Y$ be a given set. The {\em preimage} or {\em inverse image} of $\Theta$ under the mapping $F$ is defined by
\begin{equation*}
	F^{-1}(\Theta)=\big\{x\in X\; \big |\; F(x)\cap \Theta\ne\emptyset\big\}.
\end{equation*}

The  theorem below gives us  a formula to compute the normal cone to $F^{-1}(\Theta)$  at a point of interest via the normal cone to $\Theta$ and the coderivative of $F$  at certain points.

\begin{Theorem}\label{Theo-code-preimage}
	Let $F\colon X\tto Y$ be a set-valued mapping and let $\Theta\subset Y$. Suppose that $F$ is a polyhedral convex set-valued mapping and $\Theta$ is a GPCS, or $F$ is a generalized polyhedral convex set-valued mapping and $\Theta$ is a PCS. Then for any $\bar{x}\in F^{-1}(\Theta)$ and $\bar{y}\in F(\bar{x})\cap \Theta$ we have the representation
	\begin{equation}\label{Con-inverse_a}
		N(\bar{x};F^{-1}(\Theta))=D^*F(\bar{x},\bar{y})\big(N(\bar{y};\Theta)\big).
	\end{equation}
\end{Theorem}
{\bf Proof.}  Similarly as in Example~\ref{indicator mapping}, consider the indicator mappings $\Delta_\Theta\colon Y\tto  Y$ and $\Delta_{F^{-1}(\Theta)}\colon X\tto  Y$. We see that
\begin{equation}\label{equation_5}
	\Delta_{F^{-1}(\Theta)}(x)=(\Delta_\Theta \circ F)(x)\ \; \mbox{\rm for all }\; x\in X.
\end{equation}
Then the representation~\eqref{Con-inverse_a} can be obtained by using Example~\ref{indicator mapping} and Theorem~\ref{scr} with $G:=\Delta_\Theta$.  Indeed, given any $\bar{x}\in F^{-1}(\Theta)$ and $\bar{y}\in F(\bar{x})\cap \Theta$, we set $\oz=0\in Y$. It can be easily verified that $(\ox,\oz)\in\mbox{\rm gph}(G\circ F)$ and $\oy\in M(\ox,\oz)$, where the last set is defined by~\eqref{M}. By our assumptions, $F$ is a polyhedral convex set-valued mapping and $G$ is a generalized polyhedral convex set-valued mapping or vice versa. So, fixing any $ w^*\in Y^*$, by Theorem~\ref{scr} and formula~\eqref{coderivative_Delta} we can infer that
\begin{equation*}\begin{array}{rcl}
		D^*(G\circ F)(\ox,\oz)(w^*) & = & D^*F(\ox,\oy)\circ D^*G(\oy,\oz)( w^*)\\
		& = & D^*F(\bar{x},\bar{y})\big(N(\bar{y};\Theta)\big).
	\end{array}
	\end{equation*} Since the relation $D^*(G\circ F)(\ox,\oz)(w^*)=D^*\left(\Delta_{F^{-1}(\Theta)}\right)(\ox,\oz)(w^*)$ is valid by~\eqref{equation_5}, this together with~\eqref{coderivative_Delta} establishes~\eqref{Con-inverse_a}.
$\h$

\medskip
Next, consider a function $f\colon X\to \oR$ and define the {\em sublevel  sets}
\begin{equation*}
\mathcal{L}_\gamma:=\big\{x\in X\; \big|\; f(x)\leq \gamma\big\},\ \; \gamma\in \R.
\end{equation*}
Our goal is to establish a formula for the normal cone to the sublevel sets associated with a generalized polyhedral convex function. To continue, for $\ox\in \dom(f)$ we use the following notation
\begin{equation*}
  \lambda\odot \partial f(\bar{x}):=
  \begin{cases}
  \lambda \partial f(\bar{x})&\; \text{if}\; \lambda >0,\\
  \partial^{\infty}f(\bar{x})&\; \text{if}\; \lambda=0.
  \end{cases}
\end{equation*}
Here $\partial^{\infty}f(\bar{x})$ denotes the {\em singular subdifferential} of $f$ at $\ox$ defined by
\begin{equation*}
\partial^\infty f(\ox)=\big\{x^*\in X^*\; |\; (x^*, 0)\in N((\ox, f(\ox)); \epi(f)\big\}.
\end{equation*}
It follows directly from the definition that for the epigraphical mapping $E_f\colon X\to \R$ given by $E_f(x):=[f(x), \infty)$ for all $x\in X$ we have
\begin{equation}\label{equation_6}
D^*E_f(\ox, f(\ox))(\lambda)=\lambda \odot f(\ox)\ \;  {\rm for\ all}\ \, \ox\in \dom(f),\; \lambda\geq 0.
\end{equation}

\begin{Corollary}\label{normals_sublevel} Let $f\colon X\to \oR$ be a generalized polyhedral convex function with $\ox\in \dom(f)$ and $f(\ox)=\gamma$. Then
\begin{equation*}
N(\ox; \mathcal{L}_\gamma)=\bigcup_{\lambda\geq 0}\lambda\odot\partial f(\ox).
\end{equation*}
\end{Corollary}
{\bf Proof.} Let $\Theta:=(-\infty, \gamma]$ and $F(x):=E_f(x)$ for $x\in X$. Then one has $\mathcal{L}_\gamma=F^{-1}(\Theta)$. Since $\gph(F)=\epi(f)$, we see that $F$ is a generalized convex polyhedral set-valued mapping. In addition, $\Theta$ is a polyhedral convex set.  Observe also that $$N(f(\ox); \Theta)=N(\gamma; (-\infty, \gamma])=[0, \infty).$$ Therefore,  by Theorem~\ref{Theo-code-preimage} and~\eqref{equation_6} we have
\begin{equation*}\label{Con-inverse}\begin{array}{rcl}
	N\big(\ox; \mathcal{L}_\gamma\big)=N\left(\bar{x};F^{-1}(\Theta)\right) & = & D^*F(\bar{x},f(\ox))\big(N(\bar{y};\Theta)\big)\\
	& = & \bigcup\limits_{\lambda\geq 0} D^*F(\bar{x},f(\ox))(\lambda)\\
	& = & \bigcup\limits_{\lambda\geq 0}\lambda\odot\partial f(\ox),
	\end{array}		
	\end{equation*}
which completes the proof of the corollary. $\h$

\medskip
In a more general setting, consider $m$ functions $f_i\colon X\to \oR$ together with $\gamma_i\in \R$ for $i=1, \ldots, m$. Let $\gamma:=(\gamma_1, \ldots, \gamma_m)$ and define
\begin{equation*}
\mathcal{L}_\gamma:=\big\{x\in X\; \big|\; f_i(x)\leq \gamma_i\ \, \mbox{\rm for all }\, i=1, \ldots,  m\big\}.
\end{equation*}
Let $I:=\{1, \ldots, m\}$. Given $\ox\in \bigcap\limits_{i\in I}\dom(f_i)$, define
\begin{equation*}
I(\ox):=\big\{i\in I\;\big |\; f_i(\ox)=\gamma_i\big\}.
\end{equation*}
The next theorem extends the result in Corollary~\ref{normals_sublevel} to the case where the functions involved are continuous.

\begin{Theorem}\label{thm_sublevel_set} Let $f_i\colon X\to \oR$ for $i=1, \ldots, m$, and let $\gamma:=(\gamma_1, \ldots, \gamma_m)\in \R^m$, $m\geq 2$. Suppose that among $f_i$, $i=1, \ldots, m$, there are at least $m-1$ polyhedral convex functions and the remaining one is generalized polyhedral convex. Then for any $\ox\in \mathcal{L}_\gamma$ we have
\begin{equation}\label{GSLS}
N(\ox; \mathcal{L}_\gamma)=\left\{\sum_{i\in I(\ox)}\lambda_i\odot\partial f(\ox)\; \big|\; \lambda_i\geq 0\ \, \mbox{\rm for all }\, i\in I(\ox)\right\}+\sum_{i\notin I(\ox)}N(\ox; \dom(f_i)).
\end{equation}
\end{Theorem}
{\bf Proof.} Consider the set-valued mapping $F\colon X\tto \R^m$ given by
\begin{equation*}
F(x):=[f_1(x), \infty)\times \cdots\times [f_m(x), \infty),\ \; x\in X.
\end{equation*}
It can be shown that $F$ is a generalized polyhedral convex set-valued mapping. Consider the set $\Theta:=(-\infty, \gamma_1]\times\cdots\times (-\infty, \gamma_m]$ and observe that $\Theta$ is a polyhedral convex set in~$\mathbb R^m$. Clearly, $\mathcal{L}_\gamma=F^{-1}(\Theta)$.

Define the following subsets of $X\times \R^{m}$:
\begin{eqnarray*}
\begin{array}{ll}
&\Omega_1:=\{(x, \lambda_1, \ldots, \lambda_m)\; |\; \lambda_1\geq f_1(x)\}=\epi(f_1)\times \R^{m-1},\\
&\Omega_2:=\{(x, \lambda_1, \ldots, \lambda_m)\; |\; \lambda_2\geq f_2(x)\},\\
&\cdots\\
&\Omega_m:=\{(x, \lambda_1, \ldots, \lambda_m)\; |\; \lambda_m\geq f_m(x)\}.
\end{array}
\end{eqnarray*}
By our assumptions, among these sets there are at least $m-1$ polyhedral convex sets and the remaining one is a generalized polyhedral convex. Note also that  $\gph(F)=\bigcap\limits_{i\in I} \Omega_i$. By induction, from the last equality and Theorem~\ref{P intersection rule} we have
\begin{equation*}
N\Big(\big(\ox, f_1(\ox), \ldots, f_m(\ox)\big); \gph(F)\Big)=\sum_{i\in I} N\Big(\big(\ox, f_1(\ox), \ldots, f_m(\ox)\big); \Omega_i\Big).
\end{equation*}
Thus,  by the special construction of $\Omega_i$, $i\in I$, one has  $$(x^*, -\lambda_1, \ldots, -\lambda_m)\in N\big((\ox, f_1(\ox), \ldots, f_m(\ox)); \gph(F)\big)$$ if and only if there exists $x^*_i$ for $i=1, \ldots, m$ such that $(x^*_i, -\lambda_i)\in N\big((\ox, f_i(\ox)); \epi(f_i)\big)$  for each $i\in I$ and $x^*=\sum\limits_{i\in I}x^*_i$. It follows that
\begin{equation}\label{equation_7}
D^*F(\ox, f_1(\ox), \ldots, f_m(\ox))(\lambda_1, \ldots, \lambda_m)= \sum\limits_{i\in I} \lambda_i\odot \partial f_i(\ox),
\end{equation}
provided that $\lambda_i\geq 0$ for all $i\in I$. We also see that
\begin{equation*}
N((f_1(\ox), \ldots, f_m(\ox)); \Theta)=\big\{(\lambda_1, \ldots, \lambda_m)\;\big|\;  \lambda_i\geq 0\ \,  \forall i\in I,\ \lambda_i=0\ \, \mbox{\rm if }\, i \notin I(\ox)\big\}.
\end{equation*}
By  Theorem~\ref{Theo-code-preimage} we have
\begin{equation*}
N(\ox; \mathcal{L}_\gamma)=N\big(\ox;  F^{-1}(\Theta)\big)=D^*F(\ox, f_1(\ox), \ldots, f_m(\ox))\big(N((f_1(\ox), \ldots, f_m(\ox)); \Theta)\big).
\end{equation*}
Therefore, taking into account that $\partial^\infty f_i(\ox)=N(\ox; \dom(f_i))$  for every $i\in I$, we can obtain~\eqref{GSLS} from~\eqref{equation_7}. $\h$

The next corollary provides a simplified version of~\eqref{GSLS} in the case where $f_i$ is continuous at $\ox$ for all  $i\in I$.

\begin{Corollary}\label{nslc}  Under the assumptions of Theorem~\ref{thm_sublevel_set}, assume in addition that all  the functions $f_i$ for   $i\in I$  are continuous at $\ox\in \mathcal{L}_\gamma$. Then we have
\begin{equation}\label{GSLS1}
N(\ox; \mathcal{L}_\gamma)=\left\{\sum_{i\in I(\ox)}\lambda_i\partial f(\ox)\; \big|\; \lambda_i\geq 0\ \, \mbox{\rm for all }\, i\in I(\ox)\right\}.
\end{equation}
\end{Corollary}
{\bf Proof.}  For each $i\in I$, since $f_i$ is continuous at $\ox$, we have $\ox\in \mbox{\rm int}(\dom(f_i))$. Thus,
$$\partial^\infty f_i(\ox)=N(\ox; \dom(f_i))=\{0\}.$$
  Moreover, as $\partial f_i(\ox)\neq\emptyset$ by the continuity of $f_i$, we see that $$\lambda_i\odot\partial f_i(\ox)=\lambda_i \partial f_i(\ox)$$ whenever $\lambda_i\geq 0$. Therefore, the equality~\eqref{GSLS1} follows directly from~\eqref{GSLS}. $\h$

\section{Generalized Differentiation for Optimal Value Functions}\label{Sect_4}
\setcounter{equation}{0}

In this section we consider the {\em optimal value/marginal} function given  by
\begin{equation}\label{OVF}
	\mu(x):=\inf\big\{\ph(x,y)\;\big|\;y\in F(x)\big\}
\end{equation}
for all $x\in X$, where $F\colon X\tto Y$ is a set-valued mapping and $\ph\colon X\times Y\to \oR$ is an extended-real-valued   function. For simplicity of  the presentation, we assume that $\mu(x)>-\infty$ for all $x\in X$.

\begin{Theorem}\label{mr2} Let $\mu$ be an optimal value function  of the form~\eqref{OVF}. Suppose that $\ph$ is a polyhedral convex function and $F$ is a generalized polyhedral convex set-valued mapping, or $\ph$ is a generalized polyhedral convex function and $F$ is a polyhedral convex set-valued mapping. For any $\ox\in\dom(\mu)$, consider the solution set
	\begin{equation*}
		S(\ox):=\big\{\oy\in F(\ox)\;\big|\;\mu(\ox)=\ph(\ox,\oy)\big\}.
	\end{equation*}
	If  $S(\ox)$ is nonempty, then for any $\oy\in S(\ox)$ we have
	\begin{equation}\label{vf1}
		\partial\mu(\ox)=\bigcup_{(u,v)\in\partial\ph(\ox,\oy)}\big[u+D^*F(\ox,\oy)(v)\big].
	\end{equation}
\end{Theorem}
{\bf Proof.} The inclusion ``$\supset$" in~\eqref{vf1} holds when $\varphi$ is merely a convex function and $F$ is merely a convex set-valued mapping. Indeed, given $w^*$ from the right-hand side of~\eqref{vf1}, we can find  $(u,v)\in\partial\ph(\ox,\oy)$ such that $w^*-u\in D^*F(\ox,\oy)(v)$. Hence, $(w^*-u,-v)\in N((\ox,\oy);\gph(F))$. It follows that
$$\langle (w^*-u,-v), (x-\ox,y-\oy)\rangle \quad \forall (x,y)\in \gph(F).$$ Then we have
\begin{equation*}
	\langle w^*, x-\ox\rangle \leq \langle (u,v), (x,y)-(\ox,\oy)\rangle\leq  \varphi(x,y)-\varphi (\ox,\oy)
\end{equation*} for any $x\in X$ and $y\in F(x)$. Taking the infimum of the right-hand side of the inequality $\langle w^*, x-\ox\rangle \leq \varphi(x,y)-\varphi (\ox,\oy)$ and using the condition $\oy\in S(\ox)$ yield
$$\langle w^*, x-\ox\rangle\leq \mu(x)-\mu(\ox)\rangle\quad \forall x\in X.$$ This means that $w^*\in \partial\mu(\ox)$.

Let us verify the inclusion~``$\subset$" in~\eqref{vf1} under the assumptions that $\ph$ is a polyhedral convex function and $F$ is a generalized polyhedral convex set-valued mapping. Pick an element $\oy\in S(\ox)$ and let $w^*\in\partial \mu(\ox)$ be given arbitrarily. For any $x\in X$, we have
\begin{eqnarray*}
	\begin{array}{rcl}
		\la w^*,x-\ox\ra \le\mu(x)-\mu(\ox) & = & \mu(x)-\ph(\ox,\oy)\\
		& \le & \ph(x,y)-\ph(\ox,\oy)
	\end{array}
\end{eqnarray*}
 for all $y\in F(x)$. This implies that, whenever $(x,y)\in X\times Y$, the  next inequality holds:
\begin{equation*}
	\la w^*,x-\ox\ra +\la 0,y-\oy\ra\le \big[\ph(x,y)+\delta\big((x,y);\gph F)\big]-\big[\ph(\ox,\oy)+\delta\big((\ox,\oy);\gph F\big)\big].
\end{equation*}
Hence, considering the function $f(x,y):=\ph(x,y)+\delta((x,y);\gph F)$ for $(x, y)\in X\times Y$, we have $(w^*,0)\in\partial f(\ox,\oy)$. Letting $h(x,y):= \delta((x,y);\gph F)$ for $(x, y)\in X\times Y$, we  see that
\begin{equation*}
	\epi(h)=\gph(F)\times [0, \infty)
\end{equation*}
is a generalized polyhedral convex set. Since $\ph$ is a polyhedral convex function, by the subdifferential sum rule in Theorem~\ref{ssr} one has
\begin{equation*}
	(w^*,0)\in\partial f(\ox,\oy)=\partial\ph(\ox,\oy)+\partial h(\ox,\oy)=\partial\ph(\ox,\oy)+N\big((\ox,\oy);\gph F\big).
\end{equation*}
This shows that $(w^*,0)=(u^*_1,v^*_1)+(u^*_2,v^*_2)$ for some $(u^*_1,v^*_1)\in\partial\ph(\ox,\oy)$ and $$(u^*_2,v^*_2)\in N((\ox,\oy);\gph F).$$  It follows that $v^*_2=-v^*_1$; hence $(u^*_2,-v^*_1)\in N((\ox,\oy);\gph F)$.  Thus, we get $u^*_2\in D^*F(\ox,\oy)(v^*_1)$ and therefore
\begin{equation*}
	w^*=u^*_1+u^*_2\in u^*_1+D^*F(\ox,\oy)(v^*_1).
\end{equation*}
So, the inclusion~``$\subset$" in~\eqref{vf1} is valid.

The verification of the inclusion~``$\subset$" in~\eqref{vf1} under the assumptions that $\ph$ is a generalized polyhedral convex function and $F$ is a polyhedral convex set-valued mapping can be done in the same way. Namely, in the above notations, it suffices no note that $\epi(h)$ is a polyhedral convex set. As $\ph$ is a generalized polyhedral convex function, we can appy the subdifferential sum rule in Theorem~\ref{ssr} to get the desired result.

The proof of the theorem is completed. $\h$

\medskip
From Theorem~\ref{mr2} we get the following useful {\em chain rule} for convex compositions.
\begin{Corollary}\label{monotonecom} Let $f\colon X\to\R$ be a real-valued convex function and let $\phi\colon\R\to \oR$ be a nondecreasing convex function. Take  $\ox\in X$ and  let $\oy=f(\ox)\in \dom(\phi)$.  If $f$ is a generalized polyhedral convex function and  $\phi$ is a polyhedral convex function or vice versa, then
	\begin{equation}\label{crf}
		\partial(\phi\circ f)(\ox)=\disp\bigcup_{\lambda\in\partial\phi(\oy)}\lambda\odot\partial f(\ox).
	\end{equation}
	
 If we assume in addition that $f$ is continuous at $\ox$, then
 \begin{equation}\label{crf1}
		\partial(\phi\circ f)(\ox)=\disp\bigcup_{\lambda\in\partial\phi(\oy)}\lambda\partial f(\ox).
	\end{equation}
\end{Corollary}
{\bf Proof.}  Observe that the composition $\phi\circ f$ is a convex function.  Let $\ph(x,y):=\phi(y)$ for $(x, y)\in X\times \R$  and $F(x):=[f(x),\infty)$ for $x\in X$. Since $\ph$ is nondecreasing, one has
\begin{equation*}
	(\phi\circ f)(x)=\inf_{y\in F(x)}\phi(y)=\inf_{y\in F(x)}\ph(x, y).
\end{equation*}
Hence we can let $\phi\circ f$ play the role of the optimal value function $\mu$ in Theorem~\ref{mr2}. As $f$ is a generalized polyhedral convex function and $\phi$ is a polyhedral convex function or vice versa, by Theorem~\ref{mr2} we have
\begin{equation*}
	\partial(\phi\circ f)(\ox)=\disp\bigcup_{\lambda\in\partial\phi(\oy)}D^*F(\ox,\oy)(\lambda).
\end{equation*}
Taking again into account that $\phi$ is nondecreasing yields $\lambda\ge 0$ for every $\lambda\in\partial\phi(\oy)$. It follows that
\begin{equation*}
	\partial(\phi\circ f)(\ox)=\bigcup_{\lambda\in\partial\phi(\oy)}D^*F(\ox,\oy)(\lambda)=\disp\bigcup_{\lambda\in\partial\phi(\oy)}\lambda\odot\partial f(\ox),
\end{equation*}
which implies~\eqref{crf}.

The simplified version \eqref{crf1} under the continuity of $f$ at $\ox$ follows from the observation in the proof of Corollary \ref{nslc}. $\h$

\end{document}